\documentclass[reqno,final]{amsart}
\usepackage{natbib}  
\usepackage{fancyhdr} 
\usepackage{color} 
\usepackage[pagebackref=true]{hyperref} 
\usepackage{graphicx} 

\usepackage{amsthm,amssymb,amsfonts,amsmath}
\usepackage{amscd} 
\usepackage{mathrsfs}
\usepackage{vmargin}
\usepackage{setspace}
\usepackage{paralist}
\usepackage{color}
\usepackage[normalem]{ulem}
\usepackage{stmaryrd} 
\usepackage[refpage,noprefix,intoc]{nomencl} 
\usepackage{wrapfig}
\usepackage{caption}
\usepackage{subcaption}

\definecolor{aleacolor}{rgb}{0.16,0.59,0.78}

\hypersetup{
breaklinks,
colorlinks=true,
linkcolor=aleacolor,
urlcolor=aleacolor,
citecolor=aleacolor}

\renewcommand{\headheight}{11pt}
\renewcommand{\cite}{\citet}

\theoremstyle{plain}
\newtheorem{theorem}{Theorem}[section]                                          
\newtheorem{proposition}[theorem]{Proposition}                          
\newtheorem{lemma}[theorem]{Lemma}
\newtheorem{corollary}[theorem]{Corollary}

\theoremstyle{definition}

\theoremstyle{remark}

\makeatletter \@addtoreset{equation}{section} \makeatother
\renewcommand\theequation{\thesection.\arabic{equation}}
\renewcommand\thefigure{\thesection.\arabic{figure}}
\renewcommand\thetable{\thesection.\arabic{table}}

\providecommand{\R}{}
\providecommand{\Z}{}

\renewcommand{\R}{\mathbb{R}}
\renewcommand{\Z}{\mathbb{Z}}

\newcommand{\E}[1]{{\mathbf E}\left[#1\right]}
\newcommand{\e}{{\mathbf E}}

\newcommand{\p}[1]{{\mathbf P}\left\{#1\right\}}

\newcommand{\I}[1]{{\mathbf 1}_{[#1]}}
\newcommand{\set}[1]{\left\{ #1 \right\}}
\newcommand{\Cprob}[2]{\mathbf{P}\set{\left. #1 \; \right| \; #2}}
\newcommand{\probC}[2]{\mathbf{P}\set{#1 \; \left|  \; #2 \right. }}

\newcommand\cB{\mathcal B}
\newcommand\cC{\mathcal C}

\newcommand\cM{\mathcal M}

\newcommand{\rB}{\mathrm{B}} 
\newcommand{\rC}{\mathrm{C}}

\newcommand{\rF}{\mathrm{F}}

\newcommand{\rM}{\mathrm{M}}

\newcommand{\rT}{\mathrm{T}}

\newcommand{\pran}[1]{\left(#1\right)}

\newcommand{\convdist}{\ensuremath{\stackrel{\mathrm{d}}{\rightarrow}}}

\begin{document}

\title{A probabilistic approach to block sizes in random maps} 
\author{Louigi Addario-Berry}
\address{Department of Mathematics and Statistics, McGill University \\ 
805 Sherbrooke Street West \\
Montr\'eal, Qu\'ebec, H3A 0B9 \\
Canada}
\email{louigi.addario@mcgill.ca}
\urladdr{\url{http://problab.ca/louigi}}
\date{December 13, 2018} 

\renewcommand{\subjclassname}{%
  \textup{2010} Mathematics Subject Classification} 

\begin{abstract} 
We present a probabilistic approach to block sizes in random maps, which yields 
straightforward and singularity analysis-free proofs of some results of \cite{brw,bfss,gw}. 
The proof also yields joint convergence in distribution of the rescaled size of the $k$'th largest 
$2$-connected block in a large random map, for any fixed $k \ge 2$, to a vector of Fr\'echet-type extreme order statistics. This seems to be a new result even when $k=2$. 
\end{abstract}

\maketitle
\section{Introduction}\label{sec:intro}
The paper \cite{bfss} is reasonably called the culmination of an extended line of research into core sizes in large random planar maps. 
The paper is an analytic {\em tour de force}, proceeding via singularity analysis of generating functions and the coalescing saddlepoint method. \citet{bfss} demonstrate how this powerful set of tools can be used to derive to local limit theorems and sharp upper and lower tail estimates. In particular, their theorems unify and strengthen the results from \cite{brw} and \cite{gw}. 

The purpose of this note is to explain a probabilistic approach to the study of large blocks in large random maps. We end up proving two results. One is a weakening of \citep[Theorem 7]{bfss}, the other a strengthening of \citep[Proposition 5]{bfss}. The main point, though, is that our approach, which is to reduce the problem to a question about outdegrees in conditioned Galton-Watson trees, feels direct and probabilistically natural (and short).  A related technique for studying various observables of ``decomposable'' random combinatorial strucutres, using Boltzmann samplers, was introduced in \cite{kostas}. We discuss the relation between our approach and that of \cite{kostas} in Section~\ref{sec:random}. 

The remainder of the introduction lays out the definitions required for the remainder of the work. Section~\ref{sec:composite} recalls Tutte's compositional approach to planar map enumeration \cite{tutte}, and describes an associated tree decomposition of maps into higher connectivity submaps. 
Randomness finally arrives in Section~\ref{sec:random}, which also contains the statements and proofs of this work's proposition, corollary, and theorem.

\subsection{Notation for maps and trees} 
We refer the reader to \cite{lz} for a careful treatment of maps on surfaces, but provide all the definitions we directly require. In this work, a (plane) map $M$ is a planar graph $(v(M),e(M))$ properly embedded in the sphere $\mathbb{S}^2$, and considered up to orientation-preserving homeomorphisms of $\mathbb{S}^2$. Here $v(M)$ and $e(M)$ are the vertices and edges of $M$, respectively. 
All maps in this work are plane, and we hereafter omit this adjective. 
We also write $\overline{e}(M)$ for the set of oriented edges of map $M$. 

We say a map $M'$ is a submap of map $M$ if $M'$  may be obtained from $M$ by removal of a subset of the vertices and a subset of the edges of $M$. Any subgraph of $(v(M),e(M))$ induces a submap of $M$, and conversely any submap of $M$ is induced by a subgraph of $(v(M),e(M))$. Note that the faces of a submap need not be faces of the original map.

\begin{figure}[ht]
\begin{subfigure}[t]{0.22\linewidth}
		\centering
		\includegraphics[width=\linewidth,page=1]{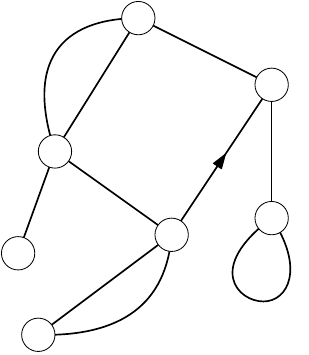}
		\caption{A map $\rM=(M,\rho)$. The root edge $\rho$ is drawn pointing from $\rho^-$ to $\rho^+$.}\label{fig:notation}		
	\end{subfigure}
	\quad
	\begin{subfigure}[t]{0.22\linewidth}
		\centering
		\includegraphics[width=\linewidth,page=2]{notation.pdf}
		\caption{The breadth-first search tree of $\rM$ has bold edges. Vertices are labelled in increasing order according to $<_{\rM}$.}\label{fig:notationtwo}
	\end{subfigure}
	\quad
	\begin{subfigure}[t]{0.22\linewidth}
		\centering
		\includegraphics[width=\linewidth,page=3]{notation.pdf}
		\caption{The oriented edges/corners are labelled in increasing order according to $\prec_{\rM}$.}
	\label{fig:notationthree}
	\end{subfigure}
	\quad
	\begin{subfigure}[t]{0.22\linewidth}
		\centering
		\includegraphics[width=\linewidth,page=6]{notation.pdf}
		\caption{The blocks of $\rM$ are shaded, and the root block has bold edges.}\label{fig:notationfour}	
	\end{subfigure}
	\caption{A map with its breadth-first search tree, corner labelling, and blocks.}
\end{figure}

A {\em rooted} map is a pair $\rM=(M,\rho)$, where $M$
is a planar map and $\rho=\rho^-\rho^+$ is an oriented edge of $M$ with tail $\rho^-$ and head $\rho^+$.  
We view $M$ as embedded in $\R^2$ so that the unbounded face lies to the right of $\rho$; this in particular gives meaning to the ``interior'' and ``exterior'' for cycles of $\rM$ (see Figure~\ref{fig:notation}). When convenient we write $v(\rM)$, etcetera, instead of $v(M)$. The {\em size} of a map is its number of edges; map $\rM$ is larger than map $\rM'$ if $|e(\rM)| \ge |e(\rM')|$. The {\em trivial} map is the map with one vertex and no edges. We root the trivial map at its unique vertex for notational convenience.

A {\em plane tree} is a connected rooted map $\rT=(T,\rho)$ with no cycles. We refer to $\rho^-$ as the root of $T$. Children and parents are then defined in the usual way. The {\em outdegree} of $v \in v(T)$ is the number of children of $v$ in $T$.

We require an ordering rule for the oriented edges of an arbitrary rooted map $\rM=(M,\rho)$. Any fixed rule would do, but for concreteness we describe a specific total order $\prec_\rM$ of $\overline{e}(M)$. Write $<_{\rM}$ for the total order of the {\em vertices} $v(M)$ induced by a breadth first search starting from $\rho^-$ using the counterclockwise order of edges around a vertex to determine exploration priority (see Figure~\ref{fig:notationtwo}). Listing the vertices according to this order as $v_1,v_2,\ldots,v_{|v(M)|}$, we in particular have $v_1=\rho^-$, $v_2=\rho^+$. We sometimes refer to $<_{\rM}$ as lexicographic order. 

Breadth-first search builds a spanning tree $\rF=\rF(\rM)$ of $\rM$ rooted at $v_1=\rho^-$: for each $v \ne \rho^-$, the parent $p(v)$ of $v$ in $\rF$ is the $<_{\rM}$-minimal neighbour $w$ of $v$. (There may be multiple edges of $\rM$ joining a node $w$ to a child $v$ of $w$, but only one of these is an edge of $\rF$; here is how to determine which. If $w=\rho^-=v_1$ then take the first copy of each edge leaving $w$ in counterclockwise order around $w$ starting from $\rho=\rho^-\rho^+$. If $w \ne \rho^-$ then take the first copy of each edge leaving $w$ in counterclockwise order starting from $wp(w)$; this makes sense inductively since $p(w) <_{\rM} w$.)

A {\em corner} of $\rM$ is a pair $(uv,uw)$ of oriented edges, where $uw$ is the successor of $uv$ in counterclockwise order around $v$. It is useful to identify oriented edges with corners: the corner corresponding to $uv$ is the corner lying to the left of its tail. This is a bijective correspondence. We define the total order $\prec_{\rM}$ on the set of corners (equivalently, the set of oriented edges) of $\rM$ as follows (see Figure~\ref{fig:notationthree}): say 
$uv \prec_{\rM} u'v'$ if either (a) $u <_{\rM} u'$ or (b) $u=u'$ and $uv$ precedes $u'v'$ in counterclockwise order around $u$ starting from $up(u)$ (or, if $u=v_1=\rho^-$, starting from $\rho$). 

\section{Planar maps as composite structures}\label{sec:composite} 

We say a rooted map $\rM$ is {\em separable} if there is a way to partition $e(\rM)$ into nonempty sets $E$ and $E'$ so that there is exactly one vertex $v$ incident to edges of both $E$ and $E'$. If $\rM$ is not separable it is called {\em $2$-connected}.\footnote{The terminology of graphs and of maps are slightly at odds here. Many graph theorists would consider the ``lollipop'' graph, with one loop and one non-loop edge, to be $2$-connected. As a map, it is not.} Write $\cM$ for the set of rooted maps, and $\cB$ for the set of $2$-connected rooted maps. \citet{tutte} showed how to count $2$-connected maps by decomposing general maps into $2$-connected submaps, then using Lagrange inversion. The remainder of the section presents this decomposition. We carefully define the tree structure associated to the decomposition, which is not explicitly used by Tutte, as it plays a key role in Section~\ref{sec:random}. 

The maximal $2$-connected submaps of $\rM$ are called the {\em blocks} of $\rM$ (hence the notation $\cB$). They are edge-disjoint, and have a natural tree structure associated to them; see Figure~\ref{fig:notationfour}. 
Write $\rB=\rB(\rM)$ for the maximal $2$-connected submap of $\rM$ containing $\rho$; call $\rB$ the {\em root block}.

\begin{wrapfigure}[19]{r}{.25\textwidth}
\rule[5pt]{0.25\textwidth}{0.5pt}
\includegraphics[width=0.25\textwidth,page=5]{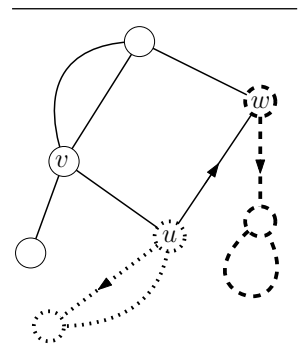}
    \caption{$\rM_{uv}$ and $\rM_{wu}$ are respectively 
    dotted and dashed.}
\label{pendants}
\rule[8pt]{0.25\textwidth}{0.5pt} 
\end{wrapfigure}

For each oriented edge $uv$ of $\rB$, there is a (possibly trivial) unique maximal submap of $\rM$ disjoint from $\rB$ except at $u$ and lying to the left of $uv$. We denote this map $\rM_{uv} =(M_{uv}, \rho_{uv})$, and call it the {\em pendant submap at $uv$} (or at the corresponding corner of $\rB$). When $\rM_{uv}$ is non-trivial, $\rho_{uv}$ is the edge of $\rM$ following $uv$ in counterclockwise order around $u$. See Figure~\ref{pendants} for an illustration. We may reconstruct $\rM$ from $\rB$ and the $2|e(\rB)|$ submaps $\{\rM_{uv}, \{u,v\} \in e(\rB)\}$ by identifying the tail of the root edge of $\rM_{uv}$ with $u \in v(\rM)$ in such a way that the root edge of $\rM_{uv}$ lies to the left of $uv$. 

\begin{figure}[ht]
\centering
\begin{subfigure}[t]{0.4\linewidth}
		\centering
		\includegraphics[width=\linewidth,page=1]{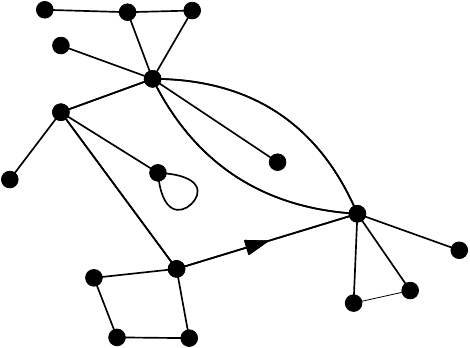}
		\caption{A map $\rM$.}\label{blocktree_edgesone}		
\end{subfigure}
	\quad
\begin{subfigure}[t]{0.4\linewidth}
		\centering
		\includegraphics[width=\linewidth,page=1]{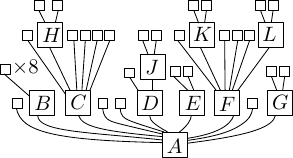}
		\caption{The tree $\rT_{\rM}$. 
Tiny squares represent trivial blocks.}\label{searchtreeone}		
\end{subfigure}
\\
\vspace{0.3cm}
\begin{subfigure}[t]{0.4\linewidth}
		\centering
		\includegraphics[width=\linewidth,page=2]{treedata.pdf}
		\caption{The decomposition of $\rM$ into blocks. Blocks are joined by grey lines according to the tree structure. Root edges of blocks are shown with arrows.}\label{blocktree_edgestwo}		
\end{subfigure}
\quad
\begin{subfigure}[t]{0.4\linewidth}
		\centering
		\includegraphics[width=\linewidth,page=1]{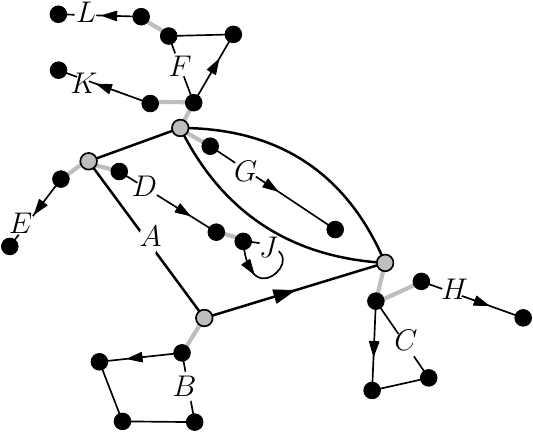}
		\caption{The correspondence between blocks and nodes of $\rT_{\rM}$. Non-trivial blocks receive the alphabetical label (from $A$ through $L$) of the corresponding node.}\label{searchtreetwo}		
\end{subfigure}
\caption{The relation between a map $\rM$ and the plane tree $\rT_{\rM}$.}
\label{blocktree_edges}
\end{figure}

Compositionally, we thereby obtain that {\em rooted maps are $2$-connected maps of rooted maps}. 
To formalize this, let $\cM_n$ (resp.\ $\cC_n$) be the set of rooted maps (resp.\ rooted $2$-connected maps) with $n$ edges, and write $M_n=|\cM_n|$, $C_n=|\cC_n|$. We take $C_0=1=M_0$. Then with $M(z) = \sum_{n \ge 0} M_n z^n$ and $C(z) = \sum_{n \ge 0} C_n z^n$, we have (see \cite{tutte}, equation (6.3))
\begin{equation}\label{eq:q_comp}
\displaywidth=\parshapelength\numexpr\prevgraf+2\relax
M(z) = C(zM(z)^2). 
\end{equation}
Now, introduce a formal variable $y$ with $y^2=z$. Then with $h(y) = y M(y^2) = z^{1/2} M(z)$, 
by (\ref{eq:q_comp}) we have $h(y) = yC(h(y)^2)$ so, by Lagrange inversion, 
\begin{equation*}
\displaywidth=\parshapelength\numexpr\prevgraf+2\relax
[z^n] M(z) = [y^{2n+1}] h(y) = \frac{1}{2n+1} [y^{2n}] C(y)^{2n+1}. 
\end{equation*}
Here is the combinatorial interpretation of this identity. Given a map $\rM=(M,\rho)$, represent the block structure of $\rM$ by the following plane tree $T_{\rM}$ defined as follows. (The construction is illustrated in Figure~\ref{blocktree_edges}.) Let $\rB=(B,\rho)$ be the block containing $\rho$, and list the {\em oriented} edges $\overline{e}(B)$ according to the order $\prec_{\rB}$ as $a_1,\ldots,a_{2|e(B)|}$. We say that the root $\emptyset$ of $T_{\rM}$ {\em represents} $\rB$ in $T_{\rM}$. 

The node $\emptyset$ has $2|e(B)|$ children in $T_{\rM}$. List them from left to right as $1,\ldots,2|e(B)|$. 
Fix $i \in \{1,\ldots,2|e(B)|\}$. If the counterclockwise successor $e_i=e_i^-e_i^+$ of $a_i$ around $a_i^-$ in $\rM$ is also in $\overline{e}(B)$ then the corner formed by $a_i$ and $e_i$ contains no pendant submap. In this case $i$ is a leaf in $T_{\rM}$. Otherwise, $e_i \in \overline{e}(M)\setminus \overline{e}(B)$. In this case write $M_i$ for the connected component of $(v(M),e(M)\setminus e(B))$ containing $\{e_i^-,e_i^+\}$, and let $\rM_i=(M_i,e_i)$. The subtree of $T_{\rM}$ rooted at $i$ is recursively defined to be the tree $T_{\rM_i}$. 
Figure~\ref{blocktree_edgesone} and~\ref{blocktree_edgestwo} show a map $\rM$ and a schematic representation of its blocktree. Figure~\ref{searchtreeone} shows the corresponding tree $\rT_{\rM}$. 

If $\rM$ is $2$-connected then $T_{\rM}$ is simply a root of outdegree $2|e(M)|$ whose children are all leaves. More generally, for each block $B$ of $\rM$, there is a corresponding node of $T_{\rM}$ with exactly $2|e(B)|$ children. In other words, given the tree $T_{\rM}$, the block sizes in $\rM$ are known. 

Given the map $\rB_{\rho}$, the map $\rM$ may be reconstructed by identifying $e_i^-$ (the tail of the root edge of $\rM_i=(M_i,e_i)$) and $a_i^-$ so that $e_i$ follows $a_i$ in counterclockwise order around $a_i^-$. (This was explained in the paragraph preceding (\ref{eq:q_comp}).)
It follows recursively that $\rM$ is uniquely specified by $T_{\rM}$ together with the set of maps $(\rB_v,v \in v(T_{\rM}))$, 
where $\rB_v$ is the block of $\rM$ represented by $v$ in $T_{\rM}$. 
If $v$ is a leaf, take $\rB_v$ to be the trivial map. Note that every node $v$ has precisely $2|e(\rB_v)|$ children in $T_{\rM}$, so $|e(T_{\rM})|=2|e(M)|$. 
For the map $\rM$ from Figure~\ref{blocktree_edgesone}, the nontrivial blocks represented by nodes of $T_{\rM}$ are shown with identifying labels in Figure~\ref{searchtreetwo}. 

\section{Random maps}\label{sec:random}
Let $\rM_n \in_u \cM_n$; this notation means that $\rM_n$ is a random variable uniformly distributed over the (finite) set $\cM_n$. 
We now describe the law of the tree $T_{\rM_n}$. Recall that $M_n = |\cM_n|$ and $C_n=|\cC_n|$, and that 
\[
M_n = \frac{2 \cdot 3^n (2n)!}{(n+2)!n!}. 
\]
Using this, the compositional equation (\ref{eq:q_comp}), and  a little thought (see \cite{tutte}, Section 6 or \cite{gj}, pages 152-153), Lagrange inversion yields 
\begin{equation}\label{eq:ck_count}
C_0=1, \quad C_k = \frac{2(3k-3)!}{k!(2k-1)!}\quad \mbox{for }k \ge 1.
\end{equation}
(The formulas for $M_n$ and $C_n$ are due to \citet{tutte}; see also \citet{brown}.) Using Stirling's approximation, the formula (\ref{eq:ck_count}) for $C_k=|\cC_k|$ implies that $C(z)$ has radius of convergence $4/27$. 
Furthermore, it is straightforward to calculate that $C(4/27)=4/3$, and that $\widehat{C}(4/27):=\sum_{k \ge 0} k(4/27)^k \cdot C_k = 4/9$. The fact $C(4/27)$ is finite is used straightaway; the second identity is noted for later use. 

Fix $z \in (0,4/27]$ and define a probability measure $\mu^z$ on the non-negative integers by 
\[
\mu^z(\{2k\}) = \frac{C_k z^k}{C(z)}. 
\]
Let $T^z$ be a Galton-Watson tree with offspring distribution $\mu^z$, and let $T^z_n$ be a random tree whose law is that of $T^z$ conditional on $|e(T^z)|=2n$.
\begin{proposition}\label{samelaw}
For all $z \in (0,4/27]$, the trees $T^z_n$ and $T_{\rM_n}$ have the same law. 
\end{proposition}
\begin{proof}[Proof of Proposition~\ref{samelaw}]
Fix a rooted plane tree $t$ with $2n$ edges, and list the outdegrees in $t$ in lexicographic order as $d_1,\ldots,d_{2n+1}$; we assume all these are even.  
We saw in Section~\ref{sec:composite} that a map $\rM$ is uniquely specified by the tree $T_{\rM}$ together with $2$-connected maps $(\rB_i,1 \le i \le 2n)$, where $\rB_i$ has $d_i/2$ edges. It follows that 
the number of maps $\rM$ with $T_{\rM}=t$ is precisely 
\[
m(t) = \prod_{i=1}^{2n} C_{\frac{d_i}{2}}. 
\]
Therefore, $\p{T_{\rM_n}=t}$ is proportional to $m(t)$. It is easily seen that this is also true for $\p{T^z_n=t}$ whatever the value of $z \in (0,4/27]$. 
\end{proof}
For the remainder of the section, let $(X_i,i \ge 1)$ be iid with law $\mu$, and write $S_k = \sum_{i=1}^k X_i$. 
Now write $\mu=\mu^{4/27}$ and $T_n = T^{4/27}_n$. 
\begin{corollary}\label{deglaw}
List the outdegrees in $T_n$ as in lexicographic order as $D_1,\ldots,D_{2n+1}$, and let $\sigma$ be a uniformly random cyclic shift of $\{1,\ldots,2n+1\}$. Then the conditional law of $(X_1,\ldots,X_{2n+1})$ given that $S_{2n+1}=2n$ is precisely that of $(D_{\sigma(1)},\ldots,D_{\sigma(2n+1)})$. 
\end{corollary}
\begin{proof}
This follows immediately from Proposition~\ref{samelaw} and the cycle lemma \citep[Lemma 6.1]{pitman}. 
\end{proof}
The corollary allows statistics about block sizes in $\rM_n$ to be deduced by studying a sequence of iid random variables conditioned on its sum. \citet{pitman} explains a quite general link between probabilistic analysis of composite structures and randomly stopped sums; he calls this {\em Kolchin's representation of Gibbs partitions}. In a sense, the point of this note is to place the study of block sizes in maps within the latter framework. 

We now state our main and only theorem. Let $A$ be a Stable$(3/2)$ random variable, characterized by its Laplace transform:
\[
\E{e^{-tA}} = e^{\Gamma(-3/2) t^{3/2}} = e^{(4\pi^{1/2}/3)t^{3/2}}.
\]
This distribution is also called a {\em map-Airy} distribution. The above scaling is used for the map-Airy distribution in (cite banderier et al); a similar scaling is used for general stable laws in (cite feller vol 2 pages 581-583). 
Also, let $(G_k,k \ge 1)$ be the ordered atoms of a rate one Poisson point process on $[0,\infty)$, so $G_k$ is $\Gamma(k)$-distributed. 
\begin{theorem}\label{main}
Let $\rM_n \in_u \cM_n$, and for $k \ge 1$ let $L_{n,k}$ be the number of edges in the $k$'th largest block of $\rM_n$. Then as $n \to \infty$, 
\[
\frac{n/3-L_{n,1}}{2^{7/6}/(27\pi)^{1/2} n^{2/3}} \convdist A, 
\]
and, jointly with the previous convergence, for any fixed $k$, as $n \to \infty$, 
\[
\quad\pran{ \frac{L_{n,j}}{(2/3)^{5/3}\pi^{-1/3}n^{2/3}},2 \le j \le k} \convdist \pran{G_{j-1}^{-3/2},2 \le j \le k}. 
\]
 \end{theorem}
Before proving the theorem, we introduce a small amount of notation. Given sequences $(Y_n)$ and $(Z_n)$ of random vectors, write $Y_n \stackrel{\mathrm{d}}{\approx} Z_n$ if $d_{\mathrm{TV}}(Y_n,Z_n) \to 0$ as $n \to \infty$, where $d_{\mathrm{TV}}$ is total variation distance. Also, for a random vector $Y$ and an event $E$, we write $(Y | E)$ for a random vector whose law is the conditional law of $Y$ given that $E$ occurs. 

\begin{proof}[Proof of Theorem~\ref{main}]
We begin with some straightforward facts about the the random variables $(X_i,i \ge 1$). 
The values of $C(4/27)$ and $\widehat{C}(4/27)$ imply that $\e{X_1}=\sum_{j\ge 0} 2j\mu(\{2j\}) = 2/3$. Furthermore, as $j \to \infty$, by Stirling's formula we have 
\[
\mu(\{2j\}) \sim \pran{\frac{8}{27\pi}}^{1/2} j^{-5/2}. 
\]
Writing $c=(\frac{8}{27\pi})^{1/2}$, it thus follows from \citep[Theorem XVII.5.2]{feller} that as $m \to \infty$, 
\begin{equation}\label{eq:airy}
\frac{S_m - 2m/3}{c m^{2/3}} \convdist A.
\end{equation}

Next, for $m \ge 1$ let $X^{m,1},\ldots,X^{m,m}$ be the decreasing rearrangement of $X_1,\ldots,X_{m}$. Then by classic results in extreme value theory (see, e.g., \citep[Section 2.2]{llr}), or by a straightforward computation, it follows that for any fixed $k$, 
\begin{equation}\label{eq:frechet}
\left(\frac{3/2}{cm}\right)^{2/3} (X^{m,i},i \le k)\convdist (G_i^{-2/3},i \ge k)\, .
\end{equation}

Now list the blocks of $\rM_n$ in decreasing order of size (number of edges) as $\rC_1,\ldots,\rC_K$, breaking ties arbitrarily, so that $L_{n,k}=|e(C_k)|$. 
By Proposition~\ref{samelaw}, the sequence $(2L_{n,k},1 \le k \le K)$ has the same law as the decreasing rearrangement of non-zero outdegrees 
in $T_n$. By Corollary~\ref{deglaw}, it follows that for all $i$ and $k$ we have 
\begin{equation}\label{eq:sizes}
\p{L_{n,k}=i} = \Cprob{X^{2n+1,k}=2i}{S_{2n+1}=2n},
\end{equation}

The large values in such collections of conditioned random variables have been studied in detail by Janson \cite{janson}. Many of the results in \cite{janson} are phrased in terms of statistics of random balls-into-boxes configurations; the connection between this and outdegrees in conditioned Galton-Watson trees is made explicit in \citep[Section 8]{janson}. One of the themes running through that work is that of {\em condensation}: for heavy-tailed random variables, conditioning a sum $S_m$ to be large is often equivalent to conditioning on having a single exceptionally large summand. 
(See \cite{al,fls,k} for other instances of this phenomenon in related settings.)

In \citep[Theorem 19.34]{janson}, Janson provides several results regarding conditional distributions such as that in (\ref{eq:sizes}). Recalling the notation introduced just before the proof, the specific result from that theorem which we use is that 
\[
((X^{2n+1,1},\ldots,X^{2n+1,2n+1})~|~ S_{2n+1}=2n) \stackrel{\mathrm{d}}{\approx} \Big(2n-S_{2n}, X^{2n,1},\ldots,X^{2n,2n}\Big)\, .
\]
For expository purposes, we include a proof of this result (in Proposition~\ref{prop:janson}, below), which closely follows that in \cite{janson}. 

Together with (\ref{eq:airy}) and (\ref{eq:sizes}), the asymptotic distributional equivalence of the first coordinate above implies that 
\[
\frac{n/3-L_{n,1}}{2^{-1/3}c n^{2/3}} \stackrel{\mathrm{d}}{\approx}  \frac{S_{2n}-2(2n/3)}{c(2n)^{2/3}} \convdist A. 
\]
The first convergence follows since $2^{-1/3}c = 2^{7/6}/(27\pi)^{1/2}$. 
Similarly, using the above asymptotic distributional equivalence together with (\ref{eq:frechet}) and (\ref{eq:sizes}) yields that for any fixed $k\ge 2$, 
\[
\left(\frac{3/2}{2cn}\right)^{2/3}(L_{n,i},2 \le i \le k) 
 \stackrel{\mathrm{d}}{\approx} 
\left(\frac{3/2}{2cn}\right)^{2/3} (X_{2n}(i),1 \le i \le k-1) 
\convdist (G_i^{-2/3},1 \le i \le k-1)\, ,
\]
which completes the proof since $(3/(4c))^{2/3} = (3/2)^{5/3}\pi^{1/3}$. \end{proof}
\begin{proposition}\label{prop:janson}
As $n \to \infty$,
\[
((X^{2n+1,1},\ldots,X^{2n+1,2n+1})~|~ S_{2n+1}=2n) \stackrel{\mathrm{d}}{\approx} \Big(2n-S_{2n}, X^{2n+1,1},\ldots,X^{2n+1,2n}\Big)\, .
\]
\end{proposition}
We first state and prove an auxiliary lemma, before proving Proposition~\ref{prop:janson}. 
\begin{lemma}\label{lem:janson}
Fix a decreasing sequence $(\delta_n)$ with $\delta_n \to 0$ slowly. 
Let 
\[
E_n = \left\{S_{2n+1}=2n,|X^{2n+1,1}-2n/3| < \delta_n n,X^{2n+1,2}<n/10\right\}\, .
\]
If $\delta_n \to 0$ sufficiently slowly then $\probC{E_n}{S_{2n+1}=2n} \to 1$. 
\end{lemma} 
\begin{proof} 
Recall that the $X_i$ are iid with $\p{X_i = 2m} \sim c m^{-5/2}$ and $\e{X_i}=2/3$. 

Write $N = \#\{i \le 2n+1: X_i \ge n/10\}$. 
By symmetry, if $N=1$ then each entry of $(X_i,1 \le i \le 2n+1)$ is equally likely to be the unique maximum. Also, for $n$ large, if $|X_i-2n/3| < \delta_n n$ then $X_i \ge n/10$. 
Provided $\delta_n \to 0$ sufficiently slowly, by the law of large numbers, $\p{|S_{2n}-4n/3| < \delta_n n}~\to~1$,~so 
\begin{align*}
\p{E_n}
& = \p{S_{2n+1}=2n,|X^{2n+1,1}-2n/3| < \delta_n n,N=1} \\
& = (2n+1) \p{S_{2n+1}=2n,|X_{2n+1}-2n/3| < \delta_n n,N=1}\\
& \ge (2n+1)\p{S_{2n+1}=2n,|X_{2n+1}-2n/3|\le \delta_n n} \\
& \ge (2n+1)\p{|S_{2n}-4n/3|\le \delta_n n} \cdot \inf_{m: |m-2n/3| \le \delta_n n} \p{X_{2n+1}=m}\\
			& \ge C n^{-3/2}\, ,
\end{align*}
for an absolute constant $C > 0$. 

In view of this lower bound on $\p{E_n}$, in order to prove the lemma it suffices to establish that $\p{S_{2n+1}=2n,E_n^c}=o(n^{-3/2})$. 
We first bound the probability that $S_{2n+1}=2n$ and $N=1$ but $E_n$ does not occur: 
\begin{align*}
& \p{S_{2n+1}=2n,N=1,|X^{2n+1,1}-2n/3| \ge \delta_n n} \\
&~=~ (2n+1)\p{S_{2n+1}=2n,X^{2n+1,2} < n/10,X_{2n+1} \ge n/10,|X_{2n+1}-2n/3| \ge \delta_n n} \\
&~=~(2n+1) \sum_{m \ge n/10: |m-2n/3| \ge \delta_n n} \p{X_{2n+1}=m} \p{S_{2n}=2n-m}\, \\
& \le (2n+1) \p{|S_{2n}-4n/3|\ge \delta_n n} \sup_{m \ge n/10: |m-2n/3| \ge \delta_n n} \p{X_{2n+1}=m}\, \\
& = o(n^{-3/2})\, ,
\end{align*}
the last bound holding since $\p{|S_{2n}-4n/3|\ge \delta_n n} \to 0$. 

It remains to prove that $\p{S_{2n+1}=2n,N \ne 1} = o(n^{-3/2})$. 
The case $N \ge 2$ is simpler: since $\p{X_1 \ge m}= O(m^{-3/2})$, 
\[
\p{S_{2n+1}=2n,N \ge 2} \le\p{N \ge 2} \le {2n+1 \choose 2} \p{X_1>  n/10,X_2 \ge n/10} = O(n^{-3})\, .
\]
In order to bound $\p{S_{2n+1}=2n,N=0}$, write $X_i' = X_i \I{X_i < n/10}$, and $S' = \sum_{i \le 2n+1} X_i'$. Then for any $t > 0$, by Markov's inequality and the independence of the $X_i'$, 
\[
\p{S_{2n+1}=2n,N=0}  = \p{S' = 2n}  \le e^{-2nt} \e e^{tS'} = e^{-2nt}\cdot \left(\e e^{tX_1'}\right)^{2n+1}. 
\]
We apply this with $t=3\log n /n$. 
To bound $\e e^{tX_1'}$, we use that for $x \in [0,5]$, $e^{x}-1-x=O(x^2)$. We thus have 
\begin{align*}
\e{e^{tX_1'}} & = 1 + t\e X_1' + \sum_{k < n/10} \p{X_1=k}(e^{t(k-1)}-1-tk)\\
& \le 1+2t/3 + C\sum_{k< n/10: tk \le 5} k^{-5/2} (tk)^2 + C \sum_{k< n/10: tk > 5} k^{-5/2} e^{tk}\, .
\end{align*}
The first sum on the final line is $O(t^{3/2})=o(1/n)$. For the second note that when $tk > 5$, 
\[
\frac{k^{-5/2} e^{tk}}{(k+1)^{-5/2} e^{t(k+1)}} = \pran{1+\frac{1}{k}}^{5/2} e^{-t} < e^{5/(2k)-t} < e^{-t/2}\, , 
\]
so the second sum is bounded by 
\[
e^{tn/10}(n/10)^{-5/2} \sum_{i \ge 0} e^{-it/2} = O(n^{-5/2} e^{tn/10}/t)=o(1/n)\, .
\]
Thus, for $t=3\log n /n$ we obtain that $\e{e^{tX_1'}}= 1+2\log n/n + o(1/n)$, so 
\[
\p{S_{2n+1}=2n,N=0} \le e^{-2nt}\cdot \left(\e e^{tX_1'}\right)^{2n+1}
= n^{-6} \pran{1+\frac{2\log n+o(1)}{n}}^{2n+1}=o(n^{3/2})\, .
\]
This completes the proof.
\end{proof}

\begin{proof}[Proof of Proposition~\ref{prop:janson}]
For $1 \le i \le 2n+1$ write $E_{n,i}=E_n \cap \{X_i = X^{2n+1,1}\}$. Then let 
\[
\begin{split} 
A = &  \Bigg\{ (x_1,\ldots,x_{2n+1}) \in \Z^{2n+1}~:   \\
&\qquad \forall i \le 2n,~0 \le x_{i} < n/10,~\Big|\sum_{i=1}^{2n} x_i - 4n/3\Big| \le \delta_n n,~x_{2n+1} = 2n-\sum_{i=1}^{2n} x_i\Bigg\}\, .
\end{split}
\]
For $n$ large, if $(x_1,\ldots,x_{2n+1}) \in A$ then $2n-\sum_{i\le 2n} x_i \ge (2/3-\delta_n)n \ge n/10 \ge \max_{i \le 2n} x_i$, so $E_{n,2n+1} = \{(X_1,\ldots,X_{2n+1}) \in A\}$. 

Next, let $\hat{E}_n = \{(X_1,\ldots,X_{2n},2n-S_{2n}) \in A\}$. Note that $(X_1,\ldots,X_{2n+1}) \in A$ if and only if $(X_1,\ldots,X_{2n},2n-S_{2n}) \in A$ and $X_{2n+1}=2n-S_{2n}$. 
Also, for all vectors $(x_1, \ldots,x_{2n+1}) \in A$ we have $|x_{2n+1} - 2n/3| \le \delta_n n$, and for such values $x_{2n+1}$, 
\[
\p{X_{2n+1}=x_{2n+1}} \sim c (n/10)^{-5/2} \, ,
\]
where as before we write $c= \pran{\frac{8}{27\pi}}^{1/2}$. Thus, uniformly over $B \subseteq A$, 
\begin{align*}
& \p{(X_1,\ldots,X_{2n+1}) \in B}\\
&~= \sum_{(x_1,\ldots,x_{2n+1}) \in B} \p{(X_1,\ldots,X_{2n})=(x_1,\ldots,x_{2n})} \p{X_{2n+1}=x_{2n+1}} \\
&~= (1+o(1)) c (n/10)^{-5/2} \sum_{(x_1,\ldots,x_{2n+1}) \in B} \p{(X_1,\ldots,X_{2n})=(x_1,\ldots,x_{2n})}\\
& = (1+o(1))  c (n/10)^{-5/2} \p{(X_1,\ldots,X_{2n},2n-S_{2n}) \in B}\, .
\end{align*}
It follows that 
\begin{align*}
\Cprob{(X_1,\ldots,X_{2n+1}) \in B}{E_{n,2n+1}}
& = \Cprob{(X_1,\ldots,X_{2n+1}) \in B}{(X_1,\ldots,X_{2n+1}) \in A} \\
& = \frac{\p{(X_1,\ldots,X_{2n+1}) \in B}}{\p{(X_1,\ldots,X_{2n+1}) \in A}} \\
& = (1+o(1)) \frac{\p{(X_1,\ldots,X_{2n},2n-S_{2n}) \in B}}{\p{(X_1,\ldots,2n-S_{2n}) \in A}} \\
& = (1+o(1)) \Cprob{(X_1,\ldots,X_{2n},2n-S_{2n}) \in B}{\hat{E}_{n}},
\end{align*}
so 
\[
((X_1,\ldots,X_{2n+1})~|~ E_{n,2n+1})\stackrel{\mathrm{d}}{\approx} 
((X_1,\ldots,X_{2n},2n-S_{2n})~|~ \hat{E}_{n})\, .
\]
For $n$ large, on $\hat{E}_n$ we have $2n-S_{2n} > \max_{i \le 2n} X_i=X^{2n,1}$, so 
\[
((X^{2n+1,1},\ldots,X^{2n+1,2n+1})~|~ E_{n,2n+1})
\stackrel{\mathrm{d}}{\approx} 
((2n-S_{2n},X^{2n,1},\ldots,X^{2n,2n})~|~ \hat{E}_{n})\, .
\]
By symmetry, the distribution of the decreasing rearrangement of $X_{1},\ldots,X_{2n+1}$ does not depend on the index at which the maximum occurs, so for all $1 \le i \le 2n+1$, 
\[
((X^{2n+1,1},\ldots,X^{2n+1,2n+1})~|~ E_{n,i})\stackrel{\mathrm{d}}{\approx} 
((X^{2n+1,1},\ldots,X^{2n+1,2n+1})~|~ E_n)\, .
\]
Moreover, Lemma~\ref{lem:janson} implies that 
\[
((X^{2n+1,1},\ldots,X^{2n+1,2n+1})~|~ E_n)
\stackrel{\mathrm{d}}{\approx}
((X^{2n+1,1},\ldots,X^{2n+1,2n+1})~|~ S_{2n+1}=2n) 
\, .
\]
Finally, 
\begin{align*}
\p{\hat{E}_n} & = \p{|S_{2n} - 4n/3| \le \delta_n n, \max_{i \le 2n} X_i < n/10}\\
& \ge 1-\p{|S_{2n} - 4n/3| > \delta_n n} - 2n\p{X_1 \ge n/10} \\
& = 1-o(1),
\end{align*}
provided $\delta_n \to 0$ sufficiently slowly, using the law of large numbers to bound the first probability on the final line, and the bound $\p{X_1 \ge n/10}=O(n^{-3/2})$ for the second. Together with the three preceding asymptotic distributional identities, this yields that 
\[
((X^{2n+1,1},\ldots,X^{2n+1,2n+1})~|~ S_{2n+1}=2n)
\stackrel{\mathrm{d}}{\approx} 
(2n-S_{2n},X^{2n,1},\ldots,X^{2n,2n})\, . \qedhere
\]
\end{proof}

 \noindent {\bf Remarks}
 \begin{enumerate}
 \item The second statement -- the convergence of the random variables $L_{n,k}$ after rescaling when $k \ge 2$ -- seems to be new. The fact that $(n^{-2/3}L_{n,2},n \ge 1)$ is a tight family of random variables, or in other words that the second largest block has size $O(n^{2/3})$ in probability, is proved in 
 \cite{gw} in some cases, and in \cite{bfss} in greater generality. 
\item \citet{kostas} showed how to use compositional schemas together with Boltzmann sampling techniques to derive information about maximal node degrees and block sizes in several families of random graphs. A similar method method was later used in \cite{ka} to derive bounds on maximal and near-maximal block sizes in random planar {\em graphs}. 
The method from \cite{kostas,ka} shares aspects with our own but yields slightly different information. In particular, it does not yield results on limiting distributions (which ours does), but does yield bounds on tail probabilities (which ours does not). 
 \item The convergence of $L_{n,1}$ is related to results from \cite{brw} and \cite{gw}. A stronger, local limit theorem for $L_{n,1}$, with explicit estimates on the rate of convergence, is given in \citet[Theorem 3]{bfss}. As mentioned earlier, the initial motivation for the current work was to show how results in this direction may be straightforwardly obtained by probabilistic arguments. With a little care, the definition of the block tree may be altered to accommodate any of the compositional schemas considered in \cite{bfss}. 
 \item In view of the preceding comment, the same line of argument should yield a version of the theorem (with constants altered appropriately) corresponding to any reasonable decomposition of a map into submaps of higher connectivity. Indeed, it seems that composite structures should in general fit within the current analytic framework. (Of course, the sorts of limit theorems one may expect will depend on the combinatorics of the specific problem. As far as I am aware, the fact that the combinatorics of maps always lead to $O(n^{2/3})$ fluctuations and Airy-type limits is thus far an empirical fact rather than a provable necessity.) 
 
 As pointed out by a referee, the block tree construction may be viewed as an instantiation of the ``enriched'' trees of Labelle \cite{labelle}, and the framework of enriched trees might be a natural one to use if one wished to generalize the arguments of the current paper; perhaps this might also shed some light on the questions implicit in the preceding paragraph. 
 
 \item It seems likely that instead of using the results of \cite{janson}, one could appeal to Theorem 1 of \cite{al}, using (2.7) from\cite{al} to control $L_{n,1}$. However, the language in \cite{janson} is closer to that of the current paper. 
 \end{enumerate}

Here are two final thoughts. First, as mentioned above, the paper \cite{bfss} proves a local limit theorem for $L_{n,1}$, with explicit error bounds in the rate of convergence. It would be interesting to recover such bounds by probabilistic methods. Second, that paper also proves essentially sharp bounds for the upper and lower tail probabilities of $L_{n,1}$; see Theorems 1 and 5. Similar tail bounds should apply in the more general settings of \cite{al,janson}. This seems like a fundamental question in large deviations of functions of iid random variables. The main result of \cite{dds} seems quite pertinent, but pertains specifically to sums rather than to more general functions.

\section{Acknowledgements}
I send my thanks to two anonymous referees, and to Mark Noy, for many useful comments and corrections. 

This work was largely written while I was visiting the Isaac Newton Institute for Mathematical Sciences during the Random Geometry programme, supported by EPSRC Grant Number EP/K032208/1. I would like to thank the Newton Institute and the Simons Foundation for their hospitality and support during this time. I also thank the University of Oxford and the Leverhulme Trust for their support during parts of this work. Finally, in all stages of this work my research was supported by NSERC and by FQRNT, I thank both institutions.

\end{document}